\documentclass[preprint,11pt]{elsarticle}

\usepackage{amssymb}
\usepackage{array}
\usepackage{longtable}
\usepackage{fancyhdr}
\usepackage{amssymb}

\newproof{proof}{Proof}
\newtheorem{pro}{Property}

\newtheorem{defn}{Definition}[section]
\newtheorem{exa}{Example}

\newtheorem{thm}{Theorem}[section]

\journal{Thermal Science }

\begin{document}

\begin{frontmatter}



\title{NEW LAPLACE-TYPE INTEGRAL TRANSFORM FOR SOLVING STEADY HEAT-TRANSFER PROBLEM }


\author{Shehu Maitama*, Weidong Zhao}

\address{School of Mathematics, Shandong University, Jinan, People's Republic of China}

\begin{abstract}
The fundamental purpose of this paper is to propose a new
Laplace-type integral transform (NL-TIT) for solving steady
heat-transfer problems. The proposed integral transform is a
generalization of the Sumudu, and the Laplace transforms and its
visualization is more comfortable than the Sumudu transform, the
natural transform, and the Elzaki transform. The suggested
integral transform is used to solve the steady heat-transfer
problems, and results are compared with the results of the
existing techniques.
\end{abstract}

\fntext[*]{*Corresponding author e-mail address:
smusman12@sci.just.edu.jo (S. Maitama), wdzhao@sdu.edu.cn (W.
Zhao). }
\end{frontmatter}


\section{Introduction}

 For more than 150 years, the motivation behind integral transforms is
easy to understand. The integral transforms have a
widely-applicable spirit of converting differential operators into
multiplication operators  from its original domain into another
domain. Besides, the symbolic manipulating and solving the
equation in the new domain is easier than manipulation and
solution in the original domain of the problem
\cite{R1}-\cite{R8}. The inverse integral transforms are always
used to mapped the manipulated solution back to the original
domain to obtain the required result.

In the mathematical literature, the famous classical integral
transforms used in differential equations, analysis, theory of
functions and integral transforms are the Laplace transform
\cite{R9} which was first introduced by a French mathematician
Pierre-Simon Laplace (1747-1827), the Fourier integral transform
\cite{R10} devised by another French mathematician Joseph Fourier
(1768-1830), and the Mellin integral transform \cite{R11} which
was introduced by a Finnish mathematician Hjalmar Mellin
(1854-1933). Besides, the Laplace transform, the Fourier
transform, and the Mellin integral transforms are similar, except
in different coordinates and have many applications in science and
engineering \cite{R12}. Moreover, in mathematics there are many
Laplace-types integral transforms such as the Laplace-Carson
transform used in the railway engineering \cite{R13}, the
z-transform applied in signal processing \cite{R14}, the Sumudu
transform used in engineering and many real-life problems
\cite{R15}, the Hankel's and Weierstrass transform applied in heat
and diffusion equations \cite{R16,R17}. In addition, we have the
natural transform \cite{R18} and Yang transform \cite{R19,R20}
used in many fields of physical science and engineering.

This paper aims to further introduce a suitable Laplace-type
integral transform for solving steady heat-transfer problems. We
will prove some important theorems and properties of the suggested
integral transform and illustrated their applications. In the next
section, we begin with the definition of the proposed Laplace-type
integral transform and introduce some useful theorems of the
integral transform.

\section*{Definition and Theorems}

\begin{defn}\label{def:ST}
The new Laplace-type integral transform of the function $v(t)$ of exponential order is defined over the set of functions,\\
$A=\left\{v(t):\exists\, C,\, \xi _{1} ,\, \xi_{2}
>0,\, \, \left|v(t)\right|<Ce^{\frac{\left|t\right|}{\xi_{i} } },\,\textnormal{if}\ t\in (-1)^{i} \times \left[0,\infty
\right)\right\}$,\\by the following integral:
\begin{equation}\label{eq1}
\Theta\left[v(t)\right](s,u)=V(s,u)=u\int _{0 }^{\infty }e^{-st}
\, v(ut)\, dt=\lim_{\beta\rightarrow\infty}\int _{0 }^{\beta
}e^{\frac{-st}{u}}v(t)dt, \, \, s\,>0,\,\, u>0,
\end{equation}
where $\Theta$ is the NL-TIT operator. It converges if the limit
of the integral exists, and diverges if not.

The inverse NL-TI transform is given by:

\begin{equation}
{\Theta^{-1}}\left[V(s,u)\right]=v(t), \, \,for\,\, t\geq0.
\end{equation}

Equivalently, the complex inversion formula of the NL-TI transform
is given by:

\begin{equation}\label{eq3}
v(t)={\Theta^{-1}}\left[V(s,u)\right]=\frac{1}{2\pi i}\int
_{\beta-i\infty}^{\beta+i\infty}\frac{1}{u}e^{\frac{st}{u}}V(s,u)\,
ds,\,\,t>0,
\end{equation}

where $s$ and $u$ are the NL-TI transform variables, and $\beta$
is a real constant and the integral in eq.\,(\ref{eq3}) is
computed along $s=\beta$ in the complex plane $s=x+iy$.
\end{defn}

\begin{thm}
The sufficient condition for the existence of the new Laplace-type
integral transform. If the function $v(t)$ is piecewise continues
on every finite interval $[0,t_0]$ and satisfies
\begin{equation}
|v(t)|\leq Ce^{\beta t},
\end{equation}
for all $t\in[t_0,\infty)$, and a constant $\beta$, then ${\Theta}\left[v(t)\right](s,u)$ exists for all $\frac{s}{u}>\beta$.\\
\end{thm}

Proof.

To prove the theorem, we must first show that the improper
integral converges for $\frac{s}{u}>\beta$. Without loss of
generality, we first split the improper integral into two parts
namely:
\begin{equation}\label{eq5}
\int_{0}^{\infty}e^{-\frac{st}{u}}v(t)dt=\int_{0}^{t_0}e^{-\frac{st}{u}}v(t)dt+\int_{t_0}^{\infty}e^{-\frac{st}{u}}v(t)dt.
\end{equation}

The first integral on the right hand side of eq.\,(\ref{eq5})
exists by the first hypothesis, hence the existence of the
Laplace-type integral transform completely depends on the second
integral. Then by the second hypothesis we deduce:

\begin{equation}
\left|e^{-\frac{st}{u}}v(t)\right|\leq Ce^{-\frac{st}{u}}e^{\beta
t}=Ce^{-\frac{(s-\beta u)t}{u}}.
\end{equation}

Thus

\begin{equation}\label{eq7}
\int_{t_0}^{\infty}Ce^{-\frac{(s-\beta u)t}{u}}dt=\frac{u}{s-\beta
u}Ce^{-\frac{(s-\beta u)t_0}{u}}.
\end{equation}

Hence, eq.\,(\ref{eq7}) converges for $\beta<\frac{s}{u}$. This
implies by the comparison test for improper integrals theorem,
${\Theta}\left[v(t)\right](s,u)$ exists for
$\beta<\frac{s}{u}$.This complete the proof.\qquad\qquad $\Box$

In the next theorem, we prove the uniqueness of the NL-TI
transform.

\begin{thm} Uniqueness of the new Laplace-type integral transform.

Let $v(t)$ and $w(t)$ be continuous  functions defined for
$t\geq0$ and having NL-TI transforms of $V(s,u)$ and $W(s,u)$
respectively. If $V(s,u)=W(s,u)$, then $v(t)=w(t)$.\\
\end{thm}

Proof.

From the inverse NL-TI transform eq.\,(\ref{eq3}), we have

\begin{equation}\label{eq8}
v(t)=\frac{1}{2\pi
i}\int_{\beta-i\infty}^{\beta+i\infty}\frac{1}{u}e^{\frac{st}{u}}V(s,u)\,
ds.
\end{equation}

 Since $V(s,u)=W(s,u)$ by the second hypothesis, then replacing this in eq.\,(\ref{eq8}), we obtain

\begin{equation}
v(t)=\frac{1}{2\pi
i}\int_{\beta-i\infty}^{\beta+i\infty}\frac{1}{u}e^{\frac{st}{u}}W(s,u)\,
ds.
\end{equation}

This implies

\begin{equation}\label{eq10}
v(t)=\frac{1}{2\pi
i}\int_{\beta-i\infty}^{\beta+i\infty}\frac{1}{u}e^{\frac{st}{u}}V(s,u)\,
ds=w(t).
\end{equation}

Hence, eq.\,(\ref{eq10}) proves the uniqueness of the NL-TI
transform.\qquad\qquad $\Box$

\begin{thm} Convolution theorem of the NL-TI transform.\label{thm3}
 Let the functions $v(t)$
and $w(t)$ be in set A. If $V(s,u)$\,\,and $W(s,u)$ are the NL-TI
transforms of the functions $v(t)$ and $w(t)$ respectively, then
\begin{equation}\label{eq11}
{\Theta}\left[(v*w)(t)\right]=V(s,u)W(s,u).
\end{equation}
Where $v*w$ is the convolution of two functions $v(t)$ and $w(t)$
which is defined as:
\begin{equation}\label{eq12}
\int _{0 }^{t }v(\tau)w(t-\tau)d\tau=\int _{0 }^{t
}v(t-\tau)w(\tau)d\tau.\\
\end{equation}
\end{thm}

Proof.

 Based on eq.\,(\ref{eq1}) and eq.\,(\ref{eq12}), we
get:
\begin{eqnarray*}
{\Theta}\left[\int_{0}^{t}v(\tau)w(t-\tau)\right]&=&\int _{0
}^{\infty}e^{\frac{-st}{u}}\left(\int_{0}^{t}v(\tau)w(t-\tau)\right)d\tau.
\end{eqnarray*}
Changing the limit of integration yields:
\begin{eqnarray*}
{\Theta}\left[\int_{0}^{t}v(\tau)w(t-\tau)\right]&=&\int _{0
}^{\infty}\left(v(\tau)\int_{\tau}^{\infty}e^{\frac{-st}{u}}w(t-\tau)dt\right)d\tau.
\end{eqnarray*}
Substituting $\vartheta=t-\tau$ in the inner integral, we deduce:
\begin{eqnarray*}
&&\int_{\tau}^{\infty}e^{\frac{-st}{u}}w(t-\tau)dt=\int _{0
}^{\infty}e^{-\frac{(\vartheta+\tau)}{u}s}w(\vartheta)d\vartheta\\
&=&e^{-\frac{\tau s}{u}}\int _{0}^{\infty}e^{-\frac{\vartheta
s}{u}}w(\vartheta)d\vartheta=e^{-\frac{\tau s}{u}}W(s,u).
\end{eqnarray*}
Hence,
\begin{eqnarray}
&&{\Theta}\left[\int_{0}^{t}v(\tau)w(t-\tau)\right]=\int _{0
}^{\infty}v(\tau)e^{\frac{-s\tau}{u}}W(s,u)d\tau\nonumber\\
&=&W(s,u)\int _{0
}^{\infty}v(\tau)e^{\frac{-s\tau}{u}}d\tau=V(s,u)W(s,u).
\end{eqnarray}
This complete the proof.\qquad\qquad $\Box$

\begin{thm}\label{thm4} Derivative of the NL-TI transform.
 Suppose that ${\Theta}
\left[v(t)\right](s,u)$ exists and that $v(t)$ is differentiable
n-times on the interval $(0,\infty)$ with $n^{th}$ derivative
$v^{(n)}(t)$, then

\begin{equation}
{\Theta} \left[v'(t)\right](s,u)=\frac{s}{u} V(s,u)-v(0),
\end{equation}
\begin{equation} {\Theta}\left[v''(t)\right](s,u)=\frac{s^{2} }{u^{2} }
V(s,u)-\frac{s}{u}v(0)-v'(0),
\end{equation}
\begin{equation} {\Theta}\left[v'''(t)\right](s,u)=\frac{s^{3} }{u^{3} } V(s,u)-\frac{s^2}{u^2
}v(0)-\frac{s}{u}v'(0)-v''(0),
\end{equation}
\begin{equation}
 \vdots\nonumber
\end{equation}
\begin{equation}\label{eq17}
{\Theta}
\left[v^{(n)}(t)\right](s,u)=\frac{s^n}{u^n}V(s,u)-\sum_{k=0}^{n-1}\left(\frac{s}{u}\right)^{n-(k+1)}v^{(k)}(0).
\end{equation}

Proof.

Using Definition \ref{def:ST} of the NL-TI transform and
integration by parts, we deduce:
\begin{eqnarray}
{\Theta} \left[v'(t)\right](s,u)&=&\int _{0 }^{\infty
}e^{\frac{-st}{u}} v'(t)dt\nonumber\\&=&-v(0)+\frac{s}{u}\int _{0
}^{\infty }e^{\frac{-st}{u}} v( t)dt =-v(0)+\frac{s}{u}V(s,u).
\end{eqnarray}

\begin{eqnarray}
{\Theta}\left[v''(t)\right](s,u)&=&\frac{s}{u}{\Theta}\left[v'(t)\right](s,u)-v'(0)=\frac{s}{u}\left[-v(0)+\frac{s}{u}V(s,u)\right]-v'(0)\nonumber\\&=&-v'(0)-\frac{s}{u}v(0)+\frac{s^2}{u^2}V(s,u).
\end{eqnarray}

\begin{eqnarray}
{\Theta}\left[v'''(t)\right](s,u)&=&\frac{s}{u}{\Theta}\left[v''(t)\right](s,u)-v''(0)=\frac{s}{u}\left[-v'(0)-\frac{s}{u}v(0)+\frac{s^2}{u^2}V(s,u)\right]-v''(0)\nonumber\\&=&-v''(0)-\frac{s}{u}v'(0)-\frac{s^2}{u^2}v(0)+\frac{s^3}{u^3}V(s,u).
\end{eqnarray}
\end{thm}
Finally, eq.\,(\ref{eq17}) follows using mathematical
induction.\qquad\qquad $\Box$

In the next theorems, we prove the NL-TI transform of
Riemann-Liouville fractional derivative
$^{RL}\hskip-2.5ptD_t^{\alpha}$ \cite{R6}, and the Caputo
fractional derivative $^{C}\hskip-2.5ptD_t^{\alpha}$ \cite{R6}.

\begin{thm} New Laplace-type integral transform of Riemann-Liouville fractional derivative.\label{thm5}
 If
$\alpha>0,\,n=1+[\alpha]$ and
$v(t)$,\,$I^{n-\alpha}v(t)$,\,$\frac{d}{dt}I^{n-\alpha}v(t)$,\,$\cdots,\,\frac{d^n}{dt^n}I^{n-\alpha}v(t)$,\,$^{RL}\hskip-2.5ptD^\alpha$
$v(t)\in A$, then

\begin{equation}
\Theta\left[^{RL}\hskip-2.5ptD_t^{\alpha}
v(t)\right]=\left(\frac{s}{u}\right)^\alpha {\Theta}
\left[v(t)\right]-\sum_{k=0}^{n-1}\left(\frac{s}{u}\right)^{n-k-1}\frac{d^{k-1}}{dt^{k-1}}I^{n-\alpha}v(0+),
\end{equation}
where $I^{\alpha}$ is the Riemann-Liouville fractional integral.
\end{thm}

Proof.

Since $^{RL}\hskip-2.5ptD_t^{\alpha}
v(t)=\frac{d^n}{dt^n}I^{n-\alpha}v(t)$.\,\,Let
$g(t)=I^{n-\alpha}v(t)$,\,\,then $^{RL}\hskip-2.5ptD_t^{\alpha}
v(t)=\frac{d^n}{dt^n}v(t)$. Applying the hypothesis of theorem
(\ref{thm4}), we get
\begin{eqnarray*}
&&\Theta\left[^{RL}\hskip-2.5ptD_t^{\alpha}
v(t)\right]=\left(\frac{s}{u}\right)^n {\Theta}
\left[v(t)\right]-\sum_{k=0}^{n-1}\left(\frac{s}{u}\right)^{n-k-1}v^{(k)}(0+)\\
&=&\left(\frac{s}{u}\right)^\alpha {\Theta}
\left[v(t)\right]-\sum_{k=0}^{n-1}\left(\frac{s}{u}\right)^{n-k-1}\frac{d^{k-1}}{dt^{k-1}}I^{n-\alpha}v(0+).
\end{eqnarray*}
The proof ends.\qquad\qquad $\Box$

\begin{thm} New Laplace-type integral transform of Caputo fractional derivative.\label{thm6}
Assume $\alpha>0,\,n=1+[\alpha]$ and
\\$v(t)$,\,$\frac{d}{dt}v(t)$,\,$\frac{d^2}{dt^2}v(t)$,\,$\cdots$,\,$\frac{d^n}{dt^n}v(t)$,
$^{C}\hskip-2.5ptD^\alpha$ $v(t)\in A$, then

\begin{equation}
\Theta\left[^{C}\hskip-2.5ptD_t^{\alpha}
v(t)\right]=\left(\frac{s}{u}\right)^\alpha {\Theta}
\left[v(t)\right]-\sum_{k=0}^{n-1}\left(\frac{s}{u}\right)^{\alpha-k-1}v^{(k)}(0+),
\end{equation}
where $^{C}\hskip-2.5ptD_t^{\alpha}$ is the Caputo fractional derivative.\\
\end{thm}

Proof.

Applying the Caputo fractional derivative \cite{R6} and theorem
(\ref{thm3}), we deduce:

\begin{eqnarray*}
^{C}\hskip-2.5ptD_t^{\alpha}
v(t)&=&\frac{1}{\Gamma(n-\alpha)}\int_{0}^{t}(t-\tau)^{n-\alpha-1}v^{(n)}(\tau)d\tau\\
&=&\frac{1}{(n-\alpha)}t^{n-\alpha-1}*v^{(n)}(\tau).
\end{eqnarray*}

Finally, using the hypothesis of theorem (\ref{thm4}) yields:
\begin{eqnarray*}
&&\Theta\left[^{C}\hskip-2.5ptD_t^{\alpha}
v(t)\right]=\frac{1}{\Gamma(n-\alpha)}\Theta\left[t^{n-\alpha-1}\right]\Theta\left[v^{(n)}(\tau)\right]\\&=&\left(\frac{s}{u}\right)^\alpha
{\Theta}
\left[v(t)\right]-\sum_{k=0}^{n-1}\left(\frac{s}{u}\right)^{\alpha-k-1}v^{(k)}(0+).
\end{eqnarray*}
This complete the proof.\qquad\qquad $\Box$

\section*{Some Properties of the NL-TI transform}
 \begin{pro} Linearity property of the NL-TI transform. Let $ {\Theta}\left[v(t)\right](s,u)=V(s,u)$ and $ {\Theta}\left[w(t)\right](s,u)=W(s,u)$, then
 \begin{equation}
 {\Theta}\left[\alpha v(t)+\beta w(t)\right](s,u)=\alpha{\Theta}\left[v(t)\right](s,u)+\beta{\Theta}\left[w(t)\right](s,u),
 \end{equation}
 where $\alpha$ and $\beta$ are constants.

Proof.

Linearity property follows directly from Definition
\ref{def:ST}.\qquad\qquad $\Box$

\end{pro}

 \begin{pro} Exponential Shifting Property of the NL-TI transform. Let the function $v(t)\in A$ and $\alpha$ is an arbitrary constant, then

 \begin{equation}
 {\Theta}\left[e^{\alpha t}v(t)\right](s,u)=V(s-\alpha u).\\
 \end{equation}
\end{pro}

Proof.

Using Definition \ref{def:ST} of the NL-TI transform, we get:

\begin{eqnarray}
{\Theta}\left[v(t)\right](s,u)&=&u\int _{0
}^{\infty}e^{-st}v(ut)dt.
\end{eqnarray}

Then

\begin{eqnarray}
{\Theta}\left[e^{\alpha
t}v(t)\right](s,u)&=&\int_{0}^{\infty}e^{(\alpha
t)}e^{\frac{-st}{u}}v(t)dt=\int _{0 }^{\infty}e^{-\frac{(s-\alpha
u)}{u}}v(t)dt\nonumber\\&=&u \int_{0}^{\infty}e^{-(s-\alpha
u)t}v(ut)dt={\Theta}\left[v(t)\right](s-\alpha u)=V(s-\alpha u).
\end{eqnarray}

In particular,

\begin{equation}
 {\Theta}\left[\sin(3t)e^{-4 t}\right](s,u)={\Theta}\left[\sin(3t)\right](s-(-4u))={\Theta}\left[\sin(3t)\right](s+4u).
 \end{equation}

Based on Definition \ref{def:ST}, the NL-TI transform of
$\sin(3t)$ is given by:

\begin{equation}\label{eq28}
 {\Theta}\left[\sin(3t)\right](s,u)=\int _{0 }^{\infty}e^{\frac{-st}{u}}\sin(3t)d t=\frac{1}{2i}\left(\frac{u}{s-3iu}-\frac{u}{s+3iu}\right)=\frac{3u^2}{s^2+9u^2}.
 \end{equation}

So, replacing the variable $s$ with $(s+4u)$ in eq.\,(\ref{eq28}),
we obtain:
\begin{equation}
 {\Theta}\left[\sin(3t)\right](s+4u)=\frac{3u^2}{(s+4u)^2+9u^2}=\frac{3u^2}{s^2+8us+25u^2}.
 \end{equation}

Alternatively,

\begin{equation}
 \int_{0}^{\infty}\sin(3t)e^{-4 t}e^{\frac{-st}{u}}dt=\frac{3u^2}{s^2+8us+25u^2}.
 \end{equation}

Moreover,

\begin{equation}
{\Theta}\left[te^{\alpha t}\right](s,u)=\frac{u^2}{(s-\alpha
u)^2}= \left\{
    \begin{array}{ll}
    \frac{1}{(s-\alpha)^2},\,\,u=1,\,\,\,\textnormal{Laplace
    transform}\ [6]
 &
        \\\\
        \frac{u^2}{(1-\alpha u)^2},\,\,s=1,\,\,\,\textnormal{Yang transform}\ [19]. &
        \\\\
    \end{array}
\right.
\end{equation}
This complete the proof.\qquad\qquad $\Box$

\begin{pro} New Laplace-type transform of integral. Let ${\Theta}\left[v(t)\right]=V(s,u)$ and $v(t)\in
A$, then

 \begin{equation}
 {\Theta}\left[\int_{0}^{t}v(\zeta)d\zeta\right]=\frac{u}{s}V(s,u).
 \end{equation}
\end{pro}

Proof.

Let $w(t)=\int_{0}^{t}v(\zeta)d\zeta$, then $w'(t)=v(t)$ and
$w(0)=0$. Computing the NL-TI transform of both sides, we get:

\begin{eqnarray}
{\Theta}\left[w'(t)\right](s,u)&=&{\Theta}\left[v(t)\right](s,u)=\frac{s}{u}{\Theta}\left[w(t)\right](s,u)-w(0)=V(s,u).
\end{eqnarray}

This implies
\begin{equation}
 {\Theta}\left[\int_{0}^{t}v(\zeta)d\zeta\right]=\frac{u}{s}V(s,u).
 \end{equation}
This complete the proof.\qquad\qquad $\Box$

\begin{pro} Multiple Shift Property of the NL-TI transform. Let ${\Theta}\left[v(t)\right](s,u)=V(s,u)$ and $v(t)\in
A$, then

 \begin{equation}\label{eq35}
 {\Theta}\left[t^nv(t)\right](s,u)=(-u)^n\frac{d^n}{ds^n}V(s,u).
 \end{equation}
\end{pro}

Proof.

By Definition \ref{def:ST} of the NL-TI transform and Leibniz's
rule, we obtain:
 \begin{eqnarray}\label{eq36}
 \frac{d}{d s}V(s,u)&=&\frac{d}{d s}\int _{0
}^{\infty }e^{\frac{-st}{u}}v(t)dt=\int _{0 }^{\infty
}\frac{\partial}{\partial
s}\left(e^{\frac{-st}{u}}\right)v(t)dt\nonumber\\&=&-\frac{1}{u}\int
_{0 }^{\infty
}e^{\frac{-st}{u}}tv(t)dt={\Theta}\left[tv(t)\right](s,u)=-u\frac{d}{d
s}V(s,u).
 \end{eqnarray}
Thus, eq.\,(\ref{eq36}) above proves the theorem for $n=1$. To
generalized the theorem, we apply the induction hypothesis. Let
assume the theorem holds for $n=k$ that is
\begin{eqnarray}
\int _{0 }^{\infty
}e^{\frac{-st}{u}}t^kv(t)dt=(-u)^k\frac{d^k}{ds^k}V(s,u).
\end{eqnarray}

Then

\begin{eqnarray}
\frac{d}{d s}\int _{0 }^{\infty
}e^{\frac{-st}{u}}t^kv(t)dt=(-u)^k\frac{d^{k+1}}{ds^{k+1}}V(s,u).
\end{eqnarray}

Alternatively, using Leibniz's rule, we deduce:

\begin{eqnarray}
&&\frac{d}{d s}\int _{0 }^{\infty }e^{\frac{-st}{u}}t^kv(t)dt=\int
_{0 }^{\infty }\frac{\partial}{\partial
s}\left(e^{\frac{-st}{u}}\right)t^kv(t)dt\nonumber\\&=&-\frac{1}{u}\int
_{0 }^{\infty
}e^{\frac{-st}{u}}t^{k+1}v(t)dt=(-u)^k\frac{d^{k+1}}{ds^{k+1}}V(s,u).
\end{eqnarray}

This implies
\begin{eqnarray}\label{eq40}
\int _{0 }^{\infty
}e^{\frac{-st}{u}}t^{k+1}v(t)dt=(-u)^{k+1}\frac{d^{k+1}}{ds^{k+1}}V(s,u).
\end{eqnarray}

Since, eq.\,(\ref{eq40}) holds for $n=k+1$, then by induction
hypothesis the prove is complete.\qquad\qquad $\Box$

\section*{Applications}
In this section, we illustrate the applicability of the proposed
Laplace-type integral transform on steady heat-transfer problems
to proves its efficiency and high accuracy.

\begin{exa}
Consider the following steady heat-transfer problem:
\begin{equation}\label{eq41}
-hMv(t)=\rho \Lambda_{c_p}v'(t),
\end{equation}
subject to the initial condition
\begin{equation}
v(0)=\beta,
\end{equation}

where $h-$ is the convection heat transfer coefficient, $M-$ is
the surface area of the body, $\rho-$ is the density of the body,
$\Lambda-$ is the volume, $c_p-$ is the specific heat of the
material, and
$v(t)$ is the temperature.\\

Applying the NL-TI transform on both sides of eq.\,(\ref{eq41}),
we get:
\begin{equation}
-hMV(s,u)=\rho \Lambda_{c_p}\left(\frac{s}{u}V(s,u)-v(0)\right).
\end{equation}
\end{exa}
Substituting the given initial condition and simplifying, we get:
\begin{equation}\label{eq44}
V(s,u)=\frac{\beta u}{s+\frac{hM}{\rho \Lambda_{c_p}}u}
\end{equation}
Taking the inverse NL-TI transform of eq.\,(\ref{eq44}), we get:
\begin{equation}
v(t)=\beta e^{-\frac{hM}{\rho \Lambda_{c_p}}t}.
\end{equation}
The exact solution is in excellent agreement with the result
obtained in \cite{R5,R20}.

\begin{exa}
Consider the following steady heat-transfer problem:
\begin{equation}\label{eq46}
v_t(x,t)=2v_{xx}(x,t),\,\,\,\,0<x<5,\,\,\,\,t>0.
\end{equation}
Subject to the boundary and initial conditions
\begin{equation}
v(0,t)=0,\,\,\,v(5,t)=0,\,\,\,\,v(x,0)=10\sin(4\pi x)-5\sin(6\pi
x).
\end{equation}
Applying the NL-TI transform on both sides of eq.\,(\ref{eq46}),
we deduce:
\begin{equation}
\frac{s}{u}V(x,s,u)-v(x,0)=\frac{2d^2V(x,s,u)}{dx^2}.
\end{equation}
\end{exa}
Substituting the given initial condition and simplifying, we get:
\begin{equation}\label{eq49}
\frac{2d^2V(x,s,u)}{dx^2}-\frac{s}{u}V(x,s,u)=-10\sin(4\pi
x)+5\sin(6\pi x).
\end{equation}
The general solution of eq.\,(\ref{eq49}) can be written as:
\begin{equation}\label{eq50}
V(x,s,u)=V_h(x,s,u)+V_p(x,s,u),
\end{equation}
where $V_h(x,s,u)$ is the solution of the homogeneous part which
is given by:
\begin{equation}\label{eq51}
V_h(x,s,u)=\alpha_1e^{\sqrt{\frac{s}{u}}x}+\alpha_2e^{-\sqrt{\frac{s}{u}}x},
\end{equation}
and $V_p(x,s,u)$ is the solution of the inhomogeneous part which
is given by:
\begin{equation}\label{eq52}
V_p(x,s,u)=\alpha\sin(4\pi x)+\beta\sin(6\pi x).
\end{equation}
Applying the boundary conditions on eq.\,(\ref{eq51}), yields

\begin{equation}
\alpha_1+\alpha_2=0\Rightarrow
\alpha_1e^{\sqrt{\frac{s}{u}}}+\alpha_2e^{-\sqrt{\frac{s}{u}}}=0\Rightarrow
V_h(x,s,u)=0,
\end{equation}

since, $\alpha_1=\alpha_2=0.$\\
Using the method of undetermined coefficients on the inhomogeneous
part, we get:
\begin{equation}
V_p(x,s,u)=10\sin(4\pi x)\frac{u}{s+32\pi^2 u}-5\sin(6\pi
x)\frac{u}{s+72\pi^2 u},
\end{equation}
since, $\alpha=10\frac{u}{s+32\pi^2
u},$\,\,\,\,and\,\,\,\,$\beta=-5\frac{u}{s+72\pi^2
u}$.\\\\
Then eq.\,(\ref{eq50}) will become:
\begin{equation}\label{eq55}
V(x,s,u)=10\sin(4\pi x)\frac{u}{s+32\pi^2 u}-5\sin(6\pi
x)\frac{u}{s+72\pi^2 u}.
\end{equation}
 Taking the inverse NL-TI transform of eq.\,(\ref{eq55}),
we obtain:
\begin{equation}
v(x,t)=10e^{-32\pi^2t}\sin(4\pi x)-5e^{-72\pi^2t}\sin(6\pi x).
\end{equation}
The exact solution is the same with the result obtained in
\cite{R9}.

\begin{exa}
Consider the following fractional porous medium equation:
\begin{equation}\label{eq57}
D_t^{\alpha}v(x,t)=D_x(v(x,t)D_xv(x,t)),\,\,0<\alpha\leq1,
\end{equation}
subject to the initial condition
\begin{equation}\label{eq58}
v(x,0)=x.
\end{equation}
Applying theorem (\ref{thm6}) on eq.\,(\ref{eq57}) subject to the
initial condition, we obtain:
\begin{equation}\label{eq59}
\Theta\left[v(x,t)\right]=\frac{u}{s}x+\frac{u^{\alpha}}{s^{\alpha}}\Theta\left[D_x(v(x,t)D_xv(x,t))\right].
\end{equation}
Computing the inverse NL-TI transform on both sides of
eq.\,(\ref{eq59}), we deduce:
\begin{equation}\label{eq60}
v(x,t)=x+\Theta^{-1}\left[\frac{u^{\alpha}}{s^{\alpha}}\Theta\left[D_x(v(x,t)D_xv(x,t))\right]\right].
\end{equation}
Based on the basic idea of the homotopy analysis method (see
\cite{R6} and references therein), we have:
\begin{equation}
v(x,t)=\sum_{n=0}^{\infty}p^{n}v_n(x,t).
\end{equation}
Then eq.\,(\ref{eq60}) will become:
\begin{equation}\label{eq62}
\sum_{n=0}^{\infty}p^{n}v_n(x,t)=x+p\left(\Theta^{-1}\left[\frac{u^{\alpha}}{s^{\alpha}}\Theta\left[D_x\left(\sum_{n=0}^{\infty}p^{n}H_n\right)\right]\right]\right),
\end{equation}
where $H_n,$\,\,is the He's polynomials \cite{R6} which represent
the
nonlinear terms\,\,$v(x,t)D_xv(x,t)$.\\
Some few components of the nonlinear terms $H_n$ are computed
below:

\begin{equation}
H_0=v_0v_{0x},\,\,H_1=v_0v_{1x}+v_1v_{0x},\,\,H_2=v_{0x}v_2+v_{1x}v_1+v_{2x}v_0,\,\,\cdots\nonumber
\end{equation}

On comparing the coefficients of same powers of $p$ in
eq.\,(\ref{eq62}), the we determine the following approximations:

\begin{eqnarray*}
p^{0} \colon v_0(x,t)&=&x,\\
p^{1}\colon
v_1(x,t)&=&\Theta^{-1}\left[\frac{u^{\alpha}}{s^{\alpha}}\Theta\left[D_x\left(H_0\right)\right]\right]\\
&=&\frac{t^\alpha}{\Gamma(1+\alpha)},
\end{eqnarray*}
\begin{eqnarray*}
p^{2}\colon
v_2(x,t)&=&\Theta^{-1}\left[\frac{u^{\alpha}}{s^{\alpha}}\Theta\left[D_x\left(H_1\right)\right]\right]\\
&=&0,\\
\vdots\\
p^{n}\colon
v_n(x,t)&=&\Theta^{-1}\left[\frac{u^{\alpha}}{s^{\alpha}}\Theta\left[D_x\left(H_{n-1}\right)\right]\right]\\
&=&0,\,\, \textnormal{for}\ \,\,n\geq2.
\end{eqnarray*}
Then the solution of eq.\,(\ref{eq57})-(\ref{eq58}) is given by:
\begin{eqnarray}\label{eq63}
v(x,t)&=&x+\frac{t^\alpha}{\Gamma(1+\alpha)}.
\end{eqnarray}

The result obtained in eq.\,(\ref{eq63})  is in excellent
agreement with the result obtained in \cite{R6}. The special case
of eq.\,(\ref{eq63}) when $\alpha=1$ is given by:

\begin{equation}\label{eq64}
v(x,t)=x+t.
\end{equation}
The result of eq.\,(\ref{eq64}) is in closed agreement with the
result obtained in \cite{R6,R7}.
\end{exa}

\section*{Conclusion}

In this paper, we introduced a powerful Laplace-type integral
transform for finding a solution of steady heat-transfer problems.
The proposed Laplace-type integral transform converges to both
Yang transform, and the Laplace transforms just by changing
variables. Many interesting properties of the suggested integral
transform are discussed and successfully applied to steady
heat-transfer problems. Finally, based on the efficiency and
simplicity of the Laplace-type integral transform, we conclude
that it is a powerful mathematical tool for solving many problems
in science and engineering.

\section*{Acknowledgment}

The authors would like to thank the anonymous reviewers, managing
editor, and editor in chief for their valuable help in improving
the manuscript. This work is supported by the Natural Science
Foundation of China (Grand No. 11571206). The first author
acknowledges the financial support of China Scholarship Council
(CSC) in Shandong University  (Grand No. 2017GXZ025381).

\section*{Nomenclature}

\begin{tabbing}
\hspace{2in} \= \hspace{2in} \= \kill
 Greek symbols \\$c_p-$ is the specific heat of the material,
[Jk$g^{-1}K^{-1}$]\>\> $x,t-$space co-ordinates, [m] \\$h-$ is the
convection heat transfer coefficient, [W$m^{-2}K^{-1}$]\>\>$v(t)-$
temperature, [K]
\hspace{1in} \\ 
 \>\> $v(x,t)-$temperature, [K]
\end{tabbing}



\begin{thebibliography}{99}
%
\bibitem{R1}
Lokenath D., Bhatta, D.,
\newblock Integral Transform and Their Applications,
\newblock {\sl CRC Press, Boca Raton, Fla., USA.,} 2014.

\bibitem{R2}
Srivastava, H.M., Luo, M., Raina, R.K.,
\newblock A New Integral Transform and Its Applications,
\newblock {\it Acta Mathematica Scientia}, 35 (2015), 6, pp. 1386--1400.

\bibitem{R3}
Yang, X.J.,
\newblock New Integral Transforms for Solving a Steady Heat-Transfer Problem,
\newblock {\it Thermal Science}, 21 (2017), pp. S79--S87.

\bibitem{R4}
Yang, X.J.,
\newblock A New Integral Transform with an Application in Heat-Transfer Problem,
\newblock {\it Thermal Science}, 20 (2016), pp. S677--S681.

\bibitem{R5}
Goodwine, B.,
\newblock Engineering Differential Equations: Theory and Applications,
\newblock {\sl Springer, New York, USA,} 2010.

\bibitem{R6}
Yan, L.M.,
\newblock Modified Homotopy perturbation Method Coupled with Laplace Transform for Fractional Heat Transfer and Porous Media Equations,
\newblock {\it Thermal Science}, 17 (2013), 5, pp. 1409--1414.

\bibitem{R7}
Pamuk, S.,
\newblock Solution of the Porous Media Equation by Adomian's Decomposition Method,
\newblock {\it Physics Letters A}, 344 (2005), pp. 184--188.

\bibitem{R8}
Elzaki, T.M.,
\newblock The New Integral Transform ''Elzaki transform'',
\newblock {\sl Glob. J. of Pur. and Appl. Math.}, 7 (2011), 1, pp. 57--64.

\bibitem{R9}
Spiegel, M.R.,
\newblock Theory and Problems of Laplace Transforms,
\newblock {\sl Schaum's Outline Series, McGraw--Hill, New York, USA.,} 1965.

\bibitem{R10}
Bracewell, R.N.,
\newblock The Fourier Transform and Its Applications,
\newblock {\it McGraw-Hill, Boston, Mass, USA.}, 2000.

\bibitem{R11}
Boyadjiev, L., Luchko, Y.,
\newblock Mellin Integral Transform Approach to Analyze the Multidimensional Diffusion-Wave Equations,
\newblock {\it Chaos Solitons and Fractals,}, 102 (2017), pp. 127-134.

\bibitem{R12}
Dattoli, G., Martinelli, M.R., Ricci, P.E.,
\newblock On New Families of Integral Transforms for the Solution of Partial Differential Equations,
\newblock {\it Integral Transforms and Special Functions}, 8 (2005), pp. 661--667.

\bibitem{R13}
L\'{e}vesque, M., et al.,
\newblock Numerical Inversion of the Laplace-Carson Transform Applied
to Homogenization of Randomly Reinforced Linear Viscoelastic
Media,
\newblock {\it Comput Mech.,} 40 (2007), pp. 771--789.

\bibitem{R14}
Cui, Y.L., et al.,
\newblock Application of the Z-Transform Technique to Modeling the Linear Lumped Networks in the HIE-FDTD Method,
\newblock {\it Journal of Electromagnetic Waves and Applications,} 27 (2013), 4, pp. 529--538.

\bibitem{R15}
Watugala, G.K.,
\newblock Sumudu Transform-A New Integral Transform to Solve Differential Equations and Control Engineering
Problems,
\newblock  {\sl Math Eng in Indust.}, 6 (1998), 1, pp. 319--329.

\bibitem{R16}
Shah, P.C., Thambynayagam, R.K.M.,
\newblock Application of the Finite Hankel Transform to a Diffusion Problem Without Azimuthal Symmetry,
\newblock {\it Transport in Porous Media}, 14 (1994), 3, pp. 247--264.

\bibitem{R17}
Karunakaran, V., Venugopal, T.,
\newblock The Weierstrass Transform for a Class of Generalized Functions,
\newblock {\it Journal of Mathematical Analysis and Applications}, 220 (1998), 2, pp. 508-527.

\bibitem{R18}
Belgacem, F.B.M., Silambarasan, R.,
\newblock  Theory of Natural Transform.
\newblock{\sl Math. in Eng. Sci., and Aeros.}, 3 (2012), 1, pp. 99--124.

\bibitem{R19}
Yang, X.J.,
\newblock A New Integral Transform Operator for Solving the Heat-Diffusion Problem,
\newblock {\it Applied Mathematics Letters}, 64 (2017), pp. 193--197.

\bibitem{R20}
Yang, X.J.,
\newblock A New Integral Transform Method for Solving Steady Heat-Transfer Problem,
\newblock {\it Thermal Science}, 20 (2016), pp. S639--S642.

\end{thebibliography}
\end{document}